\documentstyle[amscd,amssymb,verbatim,epsf]{amsart}

\begin{document}
\theoremstyle{plain}
\newtheorem{Thm}{Theorem}
\newtheorem{Cor}{Corollary}
\newtheorem{Ex}{Example}
\newtheorem{Con}{Conjecture}
\newtheorem{Main}{Main Theorem}
\newtheorem{Lem}{Lemma}
\newtheorem{Prop}{Proposition}

\theoremstyle{definition}
\newtheorem{Def}{Definition}
\newtheorem{Note}{Note}

\theoremstyle{remark}
\newtheorem{notation}{Notation}
\renewcommand{\thenotation}{}

\errorcontextlines=0
\numberwithin{equation}{section}
\renewcommand{\rm}{\normalshape}%
\newcommand{\Ge}{\ensuremath{\mathbb{G}}}
\newcommand{\J}{\ensuremath{\mathbb{J}}}
\newcommand{\Ji}{\ensuremath{\mathbb{J}}}
\newcommand{\Om}{\omega}
\renewcommand{\j}{\ensuremath{J}}
\renewcommand{\L}{{\cal L}}
\newcommand{\N}{{\cal N}}
\newcommand{\M}{{\cal M}}

\newcommand{\s}{{\cal S}}
\newcommand{\<}{\langle}
\renewcommand{\>}{\rangle}
\newcommand{\ga}{\gamma}
\newcommand{\pa}{\partial}
\newcommand{\de}{\delta}
\newcommand{\al}{\alpha}
\newcommand{\eps}{\epsilon}
\newcommand{\te}{\theta}
\newcommand{\et}{\eta}
\newcommand{\la}{\lambda}
\newcommand{\be}{\beta}
\newcommand{\si}{\sigma}
\newcommand{\ove}{\overline}

\title[Marginally trapped surfaces]%
   {Marginally trapped surfaces in spaces of oriented geodesics}
\author{Brendan Guilfoyle}
\address{Brendan Guilfoyle\\
     Institute of Technology Tralee\\
     Tralee\\
     County Kerry\\
     Ireland}
\email{guilfoylebrendan@@gmail.com}
\author{Nikos Georgiou }

\address{Nikos Georgiou\\
 Department of Mathematics and Statistics  \\
University of Cyprus \\
1678 Nicosia \\
Cyprus.}            
\email{georgiou.g.nicos@@ucy.ac.cy}

\date{28th May 2013}

\begin{abstract}
We investigate the geometric properties of marginally trapped surfaces (surfaces which have null mean curvature vector) in the spaces of 
oriented geodesics of Euclidean 3-space and hyperbolic 3-space, endowed with their canonical neutral Kaehler structures. 

We prove that every rank one surface in these four manifolds is marginally trapped. In the Euclidean case 
we show that Lagrangian rotationally symmetric sections are marginally trapped and construct an 
explicit family of marginally trapped Lagrangian tori. 

In the hyperbolic case we explore the relationship between marginally trapped and Weingarten surfaces, and construct examples of 
marginally trapped surfaces with various properties.
\end{abstract}

\maketitle

\section{Introduction}

In \cite{penrose} Roger Penrose defined a {\em marginally trapped} surface in four dimensional 
Minkowski space ${\mathbb R}^4_1$ as a spacelike surface with mean curvature vector $\vec{H}$ which is null at each point. Such surfaces play a 
central role in general relativity and, in particular, in the study of black holes and spacetime singularities.

More generally, marginally trapped submanifolds can be defined as follows. Let $(\M,G)$ be a pseudo-Riemannian manifold  and 
$\Sigma$ an immersed submanifold of codimension greater than one. We say that $\Sigma$ is {\it marginally trapped} if the mean curvature vector 
$\vec{H}$ of $\Sigma$ is null, i.e. $G(\vec{H},\vec{H})$ vanishes identically. 

The study of marginally trapped submanifolds has attracted much interest recently. For example,  Chen and Dillen \cite{chendillen} have classified 
all Lagrangian  marginally trapped surfaces in Lorentzian complex space forms. Anciaux and Godoy \cite{ango} have constructed examples of marginally 
trapped submanifolds of codimension two in de Sitter space $d{\mathbb S}^n$, anti de Sitter $Ad{\mathbb S}^n$ and the Lorentzian product 
${\mathbb S}^n\times {\mathbb R}$. Palmer \cite{palmer} has given a very important geometric interpretation of minimal spacelike surfaces of 
${\mathbb R}^4_1$: they minimize the area amongst marginally trapped surfaces satisfying natural boundary data.

In this paper, we investigate marginally trapped surfaces in the spaces of oriented geodesics of Euclidean and hyperbolic 3-space, denoted
${\mathbb L}({\mathbb R}^3)$ and ${\mathbb L}({\mathbb H}^3)$, respectively, endowed with their canonical neutral K\"ahler structures.

In particular, we prove:

\vspace{0.1in}

\noindent {\bf Theorem \ref{t:provesimrewih}.}
\emph{If $\Sigma$ is a rotationally symmetric Lagrangian section in ${\mathbb L}({\mathbb R}^3)$, then it is marginally trapped.}

\vspace{0.1in}

This extends a similar result in the space of oriented geodesics of any 3-dimensional space form of constant non-zero
sectional curvature \cite{AA}. Moreover, in Proposition \ref{p:mtsections} we construct examples of Lagrangian marginally trapped tori in 
${\mathbb L}({\mathbb R}^3)$. 
In the same section we prove our second result dealing with rank one surfaces:
\vspace{0.1in}

\noindent {\bf Theorem \ref{t:provesimrewih1}.}
\emph{Every surface in ${\mathbb L}({\mathbb R}^3)$ of rank one is marginally trapped.}
\vspace{0.1in}

Here, rank one refers to the rank of the projection $\pi:{\mathbb L}({\mathbb R}^3)\rightarrow {\mathbb S}^2$ restricted to the surface - the Gauss map.

In section 3 we investigate marginally trapped surfaces in the four dimensional space ${\mathbb L}({\mathbb H}^3)$ of oriented geodesics in 
hyperbolic 3-space endowed with its canonical neutral K\"ahler structure. For the rank one case, we prove the following:

\vspace{0.1in}

\noindent {\bf Theorem \ref{t:allrankonemarg}.}
\emph{Every surface in ${\mathbb L}({\mathbb H}^3)$ of rank one is marginally trapped.}

\vspace{0.1in}

Here again, rank refers to the rank of the two projections $\pi_j:{\mathbb L}({\mathbb H}^3)\rightarrow {\mathbb S}^2$ restricted to the surface.
We go on to find necessary and sufficient conditions for a rank two Lagrangian surface to be marginally trapped:

\vspace{0.1in}

\noindent {\bf Theorem \ref{t:grapmargtra0}.}
\emph{A rank two Lagrangian surface with potential function $h$ and Lagrangian angle $\phi$ is marginally trapped iff $\partial e^{-i\phi+h}$ 
is either a real-valued or an imaginary-valued function. }

\vspace{0.1in}
For the meanings of the terms used here, see Definitions \ref{d:potfunc} and \ref{d:lagangle}.

The relationship between marginally trapped Lagrangian surfaces in ${\mathbb L}({\mathbb H}^3)$ and Weingarten surfaces 
in ${\mathbb H}^3$ is investigated in Proposition \ref{p:ahgsu}. Finally we construct explicit examples of Lagrangian marginally 
trapped surfaces in ${\mathbb L}({\mathbb H}^3)$ with various properties, including examples that are not Weingarten.

Throughout, all curves, surfaces and maps are assumed to be smooth.
 

\section{Marginally trapped surfaces in ${\mathbb{L}}({\mathbb{R}}^{3})$}

\subsection{The neutral K\"ahler structure}

The space ${\mathbb{L}}({\mathbb{R}}^{3})$ of oriented lines in Euclidean 3-space can be identified with $T{\mathbb S}^2$ and
endowed with a canonical neutral K\"ahler structure $({\mathbb J},\Omega,{\mathbb G})$  
(for further details see \cite{Alekseevsky} and \cite{klbg}). 

Take local holomorphic coordinates $\xi$  on  ${\mathbb S}^2$ by
stereographic projection from the South pole onto the plane through
the equator. These coordinates lift to canonical holomorphic coordinates
($\xi$,$\eta$) on $T{\mathbb S}^2$.

\vspace{0.1in}

\begin{Prop}\cite{klbg}\label{p:metric}
The neutral metric on $T{\mathbb S}^2$ in holomorphic coordinates is
\[
{\Bbb{G}}=\frac{2i}{(1+\xi\bar{\xi})^2}\left(
  d\eta\otimes d\bar{\xi}-d\bar{\eta}\otimes d\xi
   +\frac{2(\xi\bar{\eta}-\bar{\xi}\eta)}{1+\xi\bar{\xi}}d\xi \otimes d\bar{\xi}
\right),
\]
while the symplectic structure is
\[
\Omega=\frac{2}{(1+\xi\bar{\xi})^2}\left(
  d\eta\wedge d\bar{\xi}+d\bar{\eta}\wedge d\xi
   +\frac{2(\xi\bar{\eta}-\bar{\xi}\eta)}{1+\xi\bar{\xi}}d\xi\wedge d\bar{\xi}
\right).
\]
\end{Prop}

\vspace{0.1in}

The K\"ahler structure is invariant under the action induced on the space of oriented lines by the Euclidean group, and so geometric properties in 
Euclidean 3-space have their analogues in the neutral K\"ahler setting. 

We consider 2-dimensional submanifolds, for which the following holds

\vspace{0.1in}

\begin{Prop}
A surface $\Sigma\subset{\mathbb L}({\mathbb R}^3)$ is Lagrangian, i.e. $\Omega|_\Sigma=0$, iff there exists a surface $S\subset{\mathbb R}^3$ which is orthogonal to the oriented
lines of $\Sigma$.
\end{Prop}

\vspace{0.1in}

Given a Lagrangian surface $\Sigma$, the projection $\pi:T{\mathbb S}^2\rightarrow{\mathbb S}^2$ restricted to $\Sigma$ is the Gauss map of the 
associated surface $S\subset{\mathbb R}^3$. The {\em rank} of $\Sigma$ is the rank of this projection, which can be zero, one or two. 

Rank zero Lagrangian surfaces are the oriented normal lines to planes in ${\mathbb R}^3$, rank one surfaces are the oriented normals to non-planar flat 
surfaces, while rank two surfaces are the oriented normals to  non-flat surfaces. A rank two surface is given locally by a section of  $\pi:T{\mathbb S}^2\rightarrow {\mathbb S}^2$, that is, 
it is given by a graph $\xi\mapsto(\xi,\eta=F(\xi,\bar{\xi}))$, for some complex function $F:{\mathbb C}\rightarrow{\mathbb C}$.

\subsection{Surfaces of rank two}

Let $\Sigma$ be a Lagrangian section given by the complex function $F$. The Lagrangian condition is equivalent to the existence of
the {\it support function} $r:\Sigma\rightarrow{\mathbb {R}}$ satisfying:
\[
F={\textstyle{\frac{1}{2}}}(1+\xi\bar{\xi})^2\bar{\partial}r,
\]
where $\partial$ and $\bar\partial$ denote differentiation with respect to $\xi$ and $\bar\xi$ respectively. 

Geometrically, $r$ is the distance along the oriented normal line between the point on the surface and the closest point to the origin.

An orthonormal basis for the tangent space $T_\gamma\Sigma$ to $\Sigma$ at a point $\gamma$ is given by
\[
E_{(1)}={\mathbb R}{\rm e}\left[\alpha_1\left(\frac{\partial}{\partial \xi}+\partial F\frac{\partial}{\partial \eta}
          +\partial \bar{F}\frac{\partial}{\partial \bar{\eta}}    \right) \right]_\gamma
\quad
E_{(2)}={\mathbb R}{\rm e}\left[\alpha_2\left(\frac{\partial}{\partial \xi}+\partial F\frac{\partial}{\partial \eta}
          +\partial \bar{F}\frac{\partial}{\partial \bar{\eta}}    \right) \right]_\gamma,
\]
where
\[
\alpha_1=\frac{(1+\xi\bar{\xi})e^{-{\scriptstyle \frac{1}{2}}\phi i+{\scriptstyle \frac{1}{4}}\pi i}}
      {\sqrt{2}|\sigma|^{\scriptstyle \frac{1}{2}}}
\qquad\qquad
\alpha_2=\frac{(1+\xi\bar{\xi})e^{-{\scriptstyle \frac{1}{2}}\phi i-{\scriptstyle \frac{1}{4}}\pi i}}
      {\sqrt{2}|\sigma|^{\scriptstyle \frac{1}{2}}},
\]
and we have introduced the complex function $\sigma=-\partial\bar{F}=|\sigma|e^{i\phi}$. Note that the induced metric is Lorentz and 
non-degenerate for $|\sigma|\neq0$, and that ${\mathbb G}(E_{(i)},E_{(j)})={\mbox{diag}}(1,-1)$.

The normal bundle $N_\gamma\Sigma$ to $\Sigma$ at a point $\gamma$ consists of vectors $T\in T_\gamma T{\mathbb S}^2$ of the form
\[
T={\mathbb R}e\left[\beta\left(\frac{\partial}{\partial \xi}+\left(\bar{\partial}\bar{F}+\frac{2(\bar{\xi}F+\xi\bar{F})}{1+\xi\bar{\xi}}\right)
\frac{\partial}{\partial \eta}-\partial\bar{F}\frac{\partial}{\partial \bar{\eta}}\right)\right]_\gamma,
\] 
for $\beta\in{\mathbb C}$.
\vspace{0.1in}

\begin{Prop}\label{p:null}
A normal vector $T\in N\Sigma$ is null iff, when written in the form above,
\[
\beta^2\partial\bar{F}=\bar{\beta}^2\bar{\partial}F.
\] 
\end{Prop}
\begin{pf}
The result follows directly from the metric expression in Proposition \ref{p:metric}.
\end{pf}

\vspace{0.1in}

\begin{Prop}
The second fundamental form of the Lagrangian section $\Sigma$ is the normal vector-valued symmetric matrix:
\[
A_{(ab)}=2{\mathbb R}{\rm e}\left[\beta_{ab}\left(\frac{\partial}{\partial \xi}
       +\left(\bar{\partial}\bar{F}+\frac{2(\bar{\xi}F+\xi\bar{F})}{1+\xi\bar{\xi}}\right)\frac{\partial}{\partial \eta}
          -\partial \bar{F}\frac{\partial}{\partial \bar{\eta}}    \right) \right],
\]
for $a,b=1,2$, where
\[
\beta_{11}=-\frac{-\sigma\bar{\partial}|\sigma|+|\sigma|^2[\partial\phi-ie^{i\phi}\bar{\partial}\phi
 -2(i\bar{\xi}-\xi e^{i\phi})(1+\xi\bar{\xi})^{-1}]}{8(1+\xi\bar{\xi})^{-2}e^{i\phi}|\sigma|^3},
\]
\[
\beta_{22}=-\frac{\sigma\bar{\partial}|\sigma|+|\sigma|^2[\partial\phi+ie^{i\phi}\bar{\partial}\phi
 -2(i\bar{\xi}+\xi e^{i\phi})(1+\xi\bar{\xi})^{-1}]}{8(1+\xi\bar{\xi})^{-2}e^{i\phi}|\sigma|^3},
\]
\[
\beta_{12}=\beta_{21}=-\frac{|\sigma|\partial |\sigma|}{8(1+\xi\bar{\xi})^{-2}e^{i\phi}|\sigma|^3}.
\]
\end{Prop}
\begin{pf}
Given a moving orthonormal basis $\{E_{(1)},E_{(2)}\}$ for the tangent space $T\Sigma$, the second fundamental form is defined to be the projection of the
ambient covariant derivative:
\[
A_{(ab)}^{\;\;\;\;\;j}=\;^\perp P_k^j\;E_{(a)}^l\overline{\nabla}_l\;E_{(b)}^k.
\]
The result follows by direct computation of these quantities on the given orthonormal frame.
\end{pf}

\vspace{0.1in}

\begin{Prop}
There are no Lagrangian sections for which the second fundamental form is null.
\end{Prop}
\begin{pf}
Applying Proposition \ref{p:null} to $A_{(12)}$ we find that nullity implies that
\begin{equation}\label{e:1}
e^{-{\scriptstyle{\frac{i}{2}}}\phi}\partial|\sigma|=e^{{\scriptstyle{\frac{i}{2}}}\phi}\bar{\partial}|\sigma|,
\end{equation}
while for $A_{(11)}-A_{(22)}$ nullity implies
\begin{equation}\label{e:2}
\partial\left(\frac{e^{-{\scriptstyle{\frac{i}{2}}}\phi}}{1+\xi\bar{\xi}}\right)
=\bar{\partial}\left(\frac{e^{{\scriptstyle{\frac{i}{2}}}\phi}}{1+\xi\bar{\xi}}\right).
\end{equation}
Note that this last equation is equivalent to ${\mathbb G}(H,H)={\mathbb G}({\mbox{tr}}A,{\mbox{tr}}A)=0$, that is, the mean curvature vector 
$H={\mbox{tr}}A$ is null.

Now, for a Lagrangian section we have the following identity:
\[
\partial\left[(1+\xi\bar{\xi})^2\partial\left(\frac{\bar{\sigma}}{(1+\xi\bar{\xi})^2}\right)\right]
=\bar{\partial}\left[(1+\xi\bar{\xi})^2\bar{\partial}\left(\frac{\sigma}{(1+\xi\bar{\xi})^2}\right)\right].
\]
Substituting equations (\ref{e:1}) and (\ref{e:2}) in this identity we find that $\partial\bar{\partial}\phi=0$,
from which we conclude that $\phi=\phi_1+\bar{\phi}_1$ for some holomorphic function $\phi_1=\phi_1(\xi)$. Substituting this in equation (\ref{e:2}) and 
rearranging we get
\[
-{\textstyle{\frac{i}{2}}}\partial\phi_1e^{-i\phi_1}-({\textstyle{\frac{i}{2}}}\partial\phi_1\xi-2)e^{-i\phi_1}\bar{\xi}
={\textstyle{\frac{i}{2}}}\bar{\partial}\bar{\phi}_1e^{i\bar{\phi}_1}-(-{\textstyle{\frac{i}{2}}}\bar{\partial}\bar{\phi}_1\bar{\xi}-2)e^{i\bar{\phi}_1}\xi.
\]
Now applying $\bar{\partial}\bar{\partial}$ to this equation we find that $\bar{\partial}\bar{\partial}\bar{\partial}(e^{i\bar{\phi}_1})=0$,
from which $\phi_1=\alpha+\beta\xi+\gamma\xi^2$ for complex constants $\alpha$, $\beta$ and $\gamma$. Back substitution implies that
$\beta=0$ and $\gamma=0$ and so $\phi$ is constant.  But this is impossible by equation (\ref{e:2}).
Thus no such Lagrangian section exists.
\end{pf}

\vspace{0.1in}

\begin{Def}
A surface is {\it marginally trapped} if its mean curvature vector is null. For Lagrangian sections, this  
condition is equivalent to equation (\ref{e:2}).

A Lagrangian section is {\it rotationally symmetric}
if (after a suitable rotation) the support function is invariant under $\xi\mapsto e^{i\beta}\xi$ for all $\beta\in[0,2\pi)$.
\end{Def}

\begin{Thm}\label{t:provesimrewih}
A rotationally symmetric Lagrangian section $\Sigma$ is  marginally trapped.
\end{Thm}
\begin{pf}
A Lagrangian graph is marginally trapped iff equation (\ref{e:2}) holds. This equation is easily seen to hold for rotationally symmetric surfaces i.e. 
ones with support function $r=r(|\xi|)$. 
\end{pf}

\vspace{0.1in}

\begin{Prop}\label{p:mtsections}
Marginally trapped Lagrangian sections which are separable with respect to Gauss polar coordinates $\xi=Re^{i\theta}$ 
are either rotationally symmetric or are tori given by the support function
\[
r=\pm \frac{R}{1+R^2}L(\theta)+r_0
\]
where $L$ is an arbitrary function of $\theta$ and $r_0\in{\mathbb{R}}$. 

\end{Prop}

\begin{pf}
The marginally trapped condition for Lagrangian surfaces is 
\[
\partial\left(\frac{e^{-{\scriptstyle{\frac{i}{2}}}\phi}}{1+\xi\bar{\xi}}\right)
=\bar{\partial}\left(\frac{e^{{\scriptstyle{\frac{i}{2}}}\phi}}{1+\xi\bar{\xi}}\right).
\]
This means that there exists a real-valued function $\Psi$ such that
\[
\frac{e^{-{\scriptstyle{\frac{i}{2}}}\phi}}{1+\xi\bar{\xi}}=\bar{\partial}\Psi,
\]
or, equivalently,
\[
\sigma=-{\scriptstyle{\frac{1}{2}}}\partial\left[(1+\xi\bar{\xi})^2\partial r\right]=|\sigma|e^{i\phi}
=|\sigma|(1+\xi\bar{\xi})^2\left(\partial\Psi\right)^2,
\]
for the support function $r$.

This implies that
\[
\left(\bar{\partial}\Psi\right)^2\partial\left[(1+\xi\bar{\xi})^2\partial r\right]
   -\left({\partial}\Psi\right)^2\bar{\partial}\left[(1+\xi\bar{\xi})^2\bar{\partial}r\right]=0.
\] 
We now seek separable solutions $r=r_1(R)r_2(\theta)$, $\Psi=\Psi_1(R)\Psi_2(\theta)$ where $\xi=Re^{i\theta}$. After some computation, the 
previous equation reduces to
\begin{equation}\label{e:eq2}
\alpha_1\alpha_2+\beta_1\beta_2=\gamma_1^2-\gamma_2^2,
\end{equation}
where
\[
\alpha_1=\alpha_1(R)=\frac{R^2[R(1+R^2)\ddot{r}_1-(1-3R^2)\dot{r}_1]\dot{\Psi}_1}{[R(1+R^2)\dot{r}_1-(1-R^2)r_1]\Psi_1}
\qquad \qquad
\alpha_2=\alpha_2(\theta)=\frac{r_2\Psi'_2}{r'_2\Psi_2},
\]
\[
\beta_1=\beta_1(R)=-\frac{R(1+R^2){r}_1\dot{\Psi}_1}{[R(1+R^2)\dot{r}_1-(1-R^2)r_1]\Psi_1}
\qquad \qquad
\beta_2=\beta_2(\theta)=\frac{r''_2\Psi'_2}{r'_2\Psi_2},
\]
\[
\gamma_1=\gamma_1(R)=\frac{R\dot{\Psi}_1}{\Psi_1}\qquad \qquad
\gamma_2=\gamma_2(\theta)=\frac{\Psi'_2}{\Psi_2}.
\]
Differentiating equation (\ref{e:eq2}) with respect to $R$ and $\theta$ we get
\[
\dot{\alpha}_1\alpha'_2+\dot{\beta}_1\beta'_2=0,
\]
with solution $\alpha_1=\lambda\beta_1+C_1$ and $\alpha_2=-\lambda^{-1}\beta_2+C_2$,
for $\lambda, C_1, C_2\in{\mathbb R}$. Substituting back into equation (\ref{e:eq2}) then yields
\[
\gamma^2_1=\lambda C_2\beta_1+C_3
\qquad\qquad
\gamma^2_2=\lambda^{-1}C_1\beta_2-C_1C_2+C_3,
\]
for $C_3\in{\mathbb R}$. 

Now note that at $R=0$ we must have $\alpha_1=\beta_1=\gamma_1=0$ for the solution to be well-defined. Thus $C_1=C_3=0$ and so
$\gamma_2=0$. Thus $\Psi'_2=0$ and, unless the surface is rotationally symmetric, we conclude that
\[
R(1+R^2)\dot{r}_1-(1-R^2)r_1=0.
\] 
Solving this ODE we find
\[
r=\frac{R}{1+R^2}L(\theta)+r_0.
\]
For some arbitrary function $L$ of $\theta$. Integrating we find that the associated graph function is
\[
\eta={\textstyle{\frac{1}{2}}}\left[(1-R^2)L+i(1+R^2)L'\right]e^{i\theta}.
\]

Note that $R=0$ is a circle of oriented lines and that this can be extended smoothly through $R=0$ by 
allowing negative values for $R$. Thus the surface is a double cover of the sphere branched along two circles, which is a torus.
\end{pf}

\vspace{0.1in}

\setlength{\epsfxsize}{5.0in}
\begin{center}
{\mbox{\epsfbox{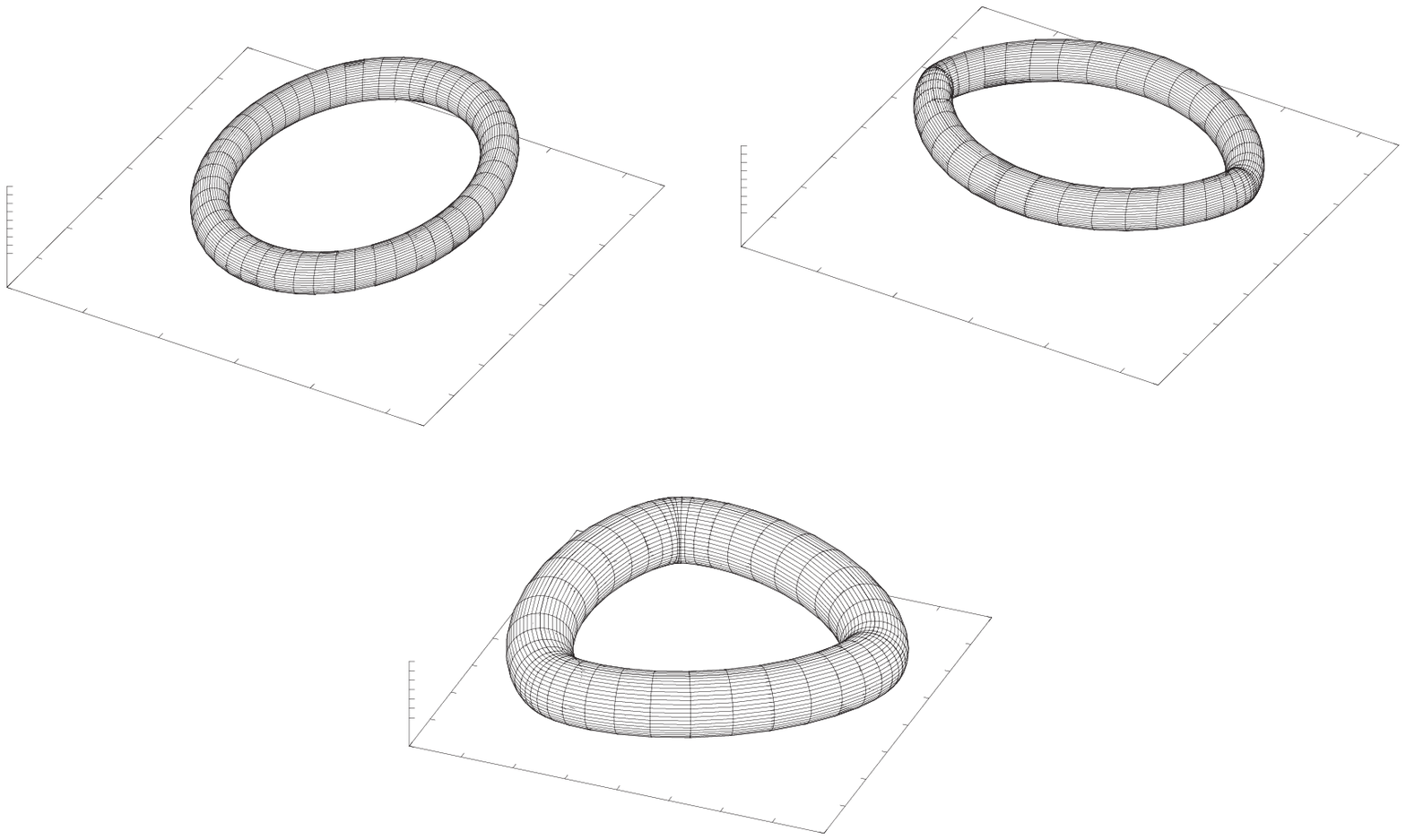}}}
\end{center}

\vspace{0.1in}

Above are pictures of the tori given by $H=a+b\cos(n\theta)$ where $a$ and $b$ are constants,  for the cases $n=0,1$ and 2. 
The resulting tori are periodic deformations of the rotationally symmetric torus obtained by squeezing it transversely.

\vspace{0.1in}

\subsection{Surfaces of rank one}

We now consider rank one surfaces in ${\mathbb L}({\mathbb R}^3)$.
\vspace{0.1in}

\begin{Thm}\label{t:provesimrewih1}
Every  surface in ${\mathbb L}({\mathbb R}^3)$ of rank one, is marginally trapped.
\end{Thm}
\begin{pf}
If $g$ is the round metric of the two-sphere ${\mathbb S}^2$, the Levi-Civita connection $\nabla$ of $g$ induces a splitting of $TT{\mathbb S}^2$ into
\begin{eqnarray}
TT{\mathbb S}^2&\simeq& T{\mathbb S}^2\oplus T{\mathbb S}^2\nonumber \\
X&\simeq& (X^h,X^v).\nonumber
\end{eqnarray}
If ${\mathbb J}, {\mathbb G}$ and $\Omega$ denote the complex structure, the neutral metric and the symplectic form respectively of the neutral K\"ahler structure endowed on $T{\mathbb S}^2$, then
\[
{\mathbb J}(X^h,X^v)=(jX^h,jX^v),\qquad \Omega((X^h,X^v),(Y^h,Y^v))=g(Y^h,X^v)-g(X^h,Y^v),
\]
where $j$ is the complex structure of $({\mathbb S}^2,g)$. Moreover, for vector fields $X,Y$, we have 
${\mathbb G}(X,Y)=\Omega({\mathbb J}X,Y)$. 

A rank one surface $\Sigma$ may be parametrised by an affine normal bundle over a regular curve $\gamma$ of ${\mathbb S}^2$ (see \cite{AGR}),
that is, they are be parametrised locally by
\begin{eqnarray}
\Psi:U\subset {\mathbb R}^2&\rightarrow& T{\mathbb S}^2\nonumber \\
(s,t)&\mapsto& (\gamma(s),V(s,t)),\nonumber
\end{eqnarray}
where $\gamma(s)$ is a regular curve of ${\mathbb S}^2$, with $s$ being the arc-length along the curve.
Therefore, $\{\gamma'(s),j\gamma'(s)\}$ is an orthonormal frame of $T{\mathbb S}^2$ along the curve $\gamma$. Write $V=a\gamma'+bj\gamma'$
and then the Frenet equations
\[
\nabla_{\gamma'}\gamma'=kj\gamma'\qquad \nabla_{\gamma'}j\gamma'=-k\gamma',
\] 
gives the derivatives of the immersion
\begin{eqnarray}
\Psi_s&=&(\gamma',\nabla_{\gamma'}V)=(\gamma',(a_s-kb)\gamma'+(b_s+ka)j\gamma')\nonumber \\
\Psi_t&=&(0,a_t\gamma'+b_t j\gamma'),\nonumber
\end{eqnarray}
where $k$ is the geodesic curvature.

If $\overline{\nabla}$ denotes the Levi-Civita connection of the ambient metric ${\mathbb G}$, a brief computation gives
\[
\overline{\nabla}_{\Psi_s}\Psi_t=(0,(a_{st}-kb_t)\gamma'+(b_{st}+ka_t)j\gamma')\quad\mbox{and}\quad \overline{\nabla}_{\Psi_t}\Psi_t=(0,0).
\]
The extrinsic curvature $h$ is defined by
\[
h(X,Y,Z)={\mathbb G}({\mathbb J}X,\overline{\nabla}_Y Z)=\Omega(X,\overline{\nabla}_Y Z),
\]
which is tri-symmetric. The only non-vanishing component is:
\begin{eqnarray}
h_{112}&=&\Omega(\Psi_s,\overline{\nabla}_{\Psi_s}\Psi_t)=-g(\gamma',(a_{st}-kb_t)\gamma'+(b_{st}+ka_t)j\gamma')\nonumber \\
&=&kb_t-a_{st}.\nonumber
\end{eqnarray}
The first fundamental form has components
\[
E=-2(b_s+ka)\qquad F=-b_t\qquad G=0,
\]
as can be seen, for example by, 
\begin{eqnarray}
E&=&{\mathbb G}(\Psi_s,\Psi_s)=\Omega({\mathbb J}\Psi_s,\Psi_s)\nonumber \\
&=&\Omega\Big((j\gamma',-(b_s+ka)\gamma'+(a_s-kb)j\gamma'),\; (\gamma',(a_s-kb)\gamma'+(b_s+ka)j\gamma')\Big)\nonumber \\
&=&-2(b_s+ka).\nonumber
\end{eqnarray}

Note that the non-degeneracy of the induced metric implies that $b_t\neq 0$. If $\vec{H}$ denotes the mean curvature vector field of the immersion $\Psi$, 
we have
\begin{eqnarray}
{\mathbb G}(2\vec{H},{\mathbb J}\Psi_s)&=&\frac{h_{111}G+h_{122}E-2h_{112}F}{EG-F^2}=\frac{2(a_{st}-kb_t)}{b_t}\nonumber \\
{\mathbb G}(2\vec{H},{\mathbb J}\Psi_t)&=&\frac{h_{112}G+h_{222}E-2h_{122}F}{EG-F^2}=0,\nonumber
\end{eqnarray}
which yields
\begin{equation}\label{e:vech}
\vec{H}=\frac{kb_t-a_{st}}{b_t}{\mathbb J}\Psi_t=\frac{kb_t-a_{st}}{b_t}(0,-b_t\gamma'+a_tj\gamma').
\end{equation}
Therefore
\begin{eqnarray}
|\vec{H}|^2&=&{\mathbb G}(\vec{H},\vec{H})\nonumber \\
&=&\left(\frac{kb_t-a_{st}}{b_t}\right)^2{\mathbb G}((0,-b_t\gamma'+a_tj\gamma'),(0,-b_t\gamma'+a_tj\gamma'))\nonumber \\
&=&\left(\frac{kb_t-a_{st}}{b_t}\right)^2\Omega((0,-a_t\gamma'-b_tj\gamma'),(0,-b_t\gamma'+a_tj\gamma'))\nonumber \\
&=&0,\nonumber 
\end{eqnarray}
and so the immersion $\Psi$ is marginally trapped.
\end{pf}

\begin{Note}
{If $a(s,t)$ is a smooth function such that
\[
a(s,t)=a_1(s)+a_2(t),
\] 
then from the expression (\ref{e:vech}), the rank one immersion $\Psi$ of the theorem above,  is minimal iff $\gamma$ is a geodesic.}
\end{Note}

\vspace{0.2in}

\section{Marginally trapped surfaces in ${\mathbb L}({\mathbb H}^3)$}

\subsection{The neutral K\"ahler structure}

In this section we summarize the geometry of the neutral K\"ahler structure $({\mathbb G},{\mathbb J},\Omega)$ on the space 
${\mathbb L}({\mathbb H}^3)$ of oriented geodesics of hyperbolic space - further details can be found in \cite{AA} and \cite{gag}.

The complex surface  $({\mathbb L}({\mathbb H}^3),{\mathbb J})$ can be identified with $({\mathbb P}^1\times {\mathbb P}^1-\overline\Delta,j\oplus j)$, 
where we have introduced the reflected diagonal, which in holomorphic coordinates $(\mu_1,\mu_2)$ on ${\mathbb P}^1\times {\mathbb P}^1$ is 
$\overline\Delta=\{(\mu_1,\mu_2)\in {\mathbb P}^1\times {\mathbb P}^1\left.{\!\!\frac{}{}}\right|\;\mu_1\bar{\mu}_2=-1\}$.

\vspace{0.1in}
\begin{Prop}
The K\"ahler metric ${\mathbb G}$ and the symplectic form $\Omega$ on  ${\mathbb L}({\mathbb H}^3)$ are expressed locally by:
\[
{\mathbb G}={\mathbb I}\mbox{m}\left[\frac{2}{(1+\mu_1\bar{\mu}_2)^2}d\mu_1\otimes d\bar{\mu}_2\right]
\qquad 
\Omega=-{\mathbb R}\mbox{e}\left[\frac{2}{(1+\mu_1\bar{\mu}_2)^2}d\mu_1\wedge d\bar{\mu}_2\right].
\]
\end{Prop}
\vspace{0.1in}

The neutral K\"ahler metric ${\mathbb G}$ is scalar flat, conformally flat and is invariant under the action induced on ${\mathbb{L}}({\mathbb{H}}^3)$ by 
the isometry group of ${\mathbb{H}}^3$. In fact, ${\mathbb G}$ is the unique K\"ahler metric on ${\mathbb{L}}({\mathbb{H}}^3)$ with this last property 
\cite{salvai}.

We consider surfaces $\Sigma\subset{\mathbb L}({\mathbb H}^3)$, where the following holds

\vspace{0.1in}

\begin{Prop}
A surface $\Sigma$ is Lagrangian, i.e. $\Omega|_\Sigma=0$, iff there exists a surface $S\subset{\mathbb H}^3$ which is orthogonal to the oriented
lines of $\Sigma$.
\end{Prop}

\vspace{0.1in}

\begin{Def}
Given an immersion $f:\Sigma\rightarrow {\mathbb L}({\mathbb H}^3)$, consider the maps $(\pi_j\circ f)_{\ast}:T\Sigma\rightarrow T{\mathbb P}^1$, for
$j=1,2$, where $\pi_j$ is projection onto the jth factor of ${\mathbb L}({\mathbb H}^3)={\mathbb{P}}^1\times {\mathbb{P}}^1-\overline{\Delta}$. 
The {\it rank} of the immersion $f$ at a point $\gamma\in\Sigma$ is defined to be the smaller of the ranks of these maps at $\gamma$, which can be 
0, 1 or 2.
\end{Def}

Note that by choice of orientation of the geodesics, the rank can be defined to be the rank of the projection onto the first factor.

A rank 0 surface in ${\mathbb L}({\mathbb H}^3)$ is the set of oriented geodesics orthogonal to a horosphere, a rank 1 surface can be locally 
parameterized by $\mu_1=\mu_1(s)$ and $\mu_2=\mu_2(s,t)$ for $(s,t)\in D\subset{\mathbb R}^2$, while  a rank 2 surface can be locally parameterized by 
$\mu_1\rightarrow (\mu_1,\mu_2=\mu_2(\mu_1,\bar{\mu}_1))$.

\subsection{Marginally trapped Weingarten surfaces}

\begin{Def}
A surface in ${\mathbb H}^3$ is said to be {\it Weingarten} if the eigenvalues of the second fundamental form (the principal curvatures) 
are functionally related. 
\end{Def}

There is an important geometric relation between Weingarten surfaces in ${\mathbb H}^3$ and the Lagrangian surfaces in ${\mathbb{L}}({\mathbb{H}}^3)$  
formed by the oriented geodesics orthogonal to $S$:

\vspace{0.1in}

\begin{Prop}\label{t:weingflat}\cite{AA} \cite{nikosbrendan}
Let $\Sigma\subset{\mathbb{L}}({\mathbb{H}}^3)$ the Lagrangian 
surface formed by the oriented geodesics orthogonal to a surface $S\subset{\mathbb{H}}^3$. 
Assume that the induced metric ${\mathbb G}_\Sigma$ is non-degenerate. 

Then $S$ is Weingarten iff the Gauss curvature of ${\mathbb G}_\Sigma$ is zero.
\end{Prop}

\vspace{0.1in}

Following the notation of Anciaux \cite{AA}, let $\phi:\tilde S\rightarrow {\mathbb H}^3$ be a Riemannian immersion and 
$g=\phi^{\ast}\<\cdot,\cdot\>$  be the induced metric. Then $S=\phi(\tilde S)$ is an immersed surface of ${\mathbb H}^3$ 
and denote the shape operator by $A$. Let $(e_1,e_2)$ be the orthonormal frame such that $Ae_i=k_i e_i$, where $k_i$ are the 
principal curvatures. The oriented normal geodesics $\Sigma$ of $S$ are the image of the map $\bar\phi:\tilde S\rightarrow {\mathbb L}({\mathbb H}^3)$ 
defined by $\bar\phi=\phi\wedge N$, where $N$ is the unit vector field along the immersed surface $S$.

In {\cite{AA}} the following is proved:

\vspace{0.1in}

\begin{Prop}\label{p:impaenri}
Assume that $S\subset{\mathbb H}^3$ either has a constant principal curvature or is a surface of revolution. Then the corresponding oriented normal 
geodesics $\Sigma$ is marginally trapped.
\end{Prop}

\vspace{0.1in}

For the converse, we prove the following:

\vspace{0.1in}

\begin{Prop}\label{p:weingmargtra}
Assume that $S$ is a Weingarten surface in ${\mathbb H}^3$ and the normal geodesics $\Sigma$ is marginally trapped. Then either $S$ has a 
constant principal curvature or it is a surface of revolution.
\end{Prop}
\begin{pf}
Let $k_1,k_2$ be the principal curvatures of $S$ and $(e_1,e_2)$ be the principal directions, that is, $Ae_i=k_ie_i$ for $i=1,2$.

From \cite{AA}, the mean curvature vector $\vec{H}$ of $\Sigma$ is 
\[
\vec{H}=(2(k_2-k_1)^2)^{-1}(e_1(k_2){\mathbb J}d\bar\phi(e_1)+e_2(k_1){\mathbb J}d\bar\phi(e_2)),
\] 
where, as before, ${\mathbb J}$ denotes the complex structure defined in ${\mathbb L}({\mathbb H}^3)$. If $\bar g=\bar\phi^{\ast}{\mathbb G}$ denotes the induced metric then
\[
\bar g (e_1, e_1)=\bar g (e_2, e_2)=0
\qquad\qquad
\bar g (e_1, e_2)=k_2-k_1.
\]
Thus,
\[
|\vec{H}|^2=\frac{e_1(k_2)e_2(k_1)}{2(k_2-k_1)}=0,
\]
which gives
\begin{equation}\label{e:fistequalll}
e_1(k_2)e_2(k_1)=0.
\end{equation}
Also the fact that $S$ is Weingarten, yields 
\begin{equation}\label{e:secequalll}
e_1(k_1)e_2(k_2)=0.
\end{equation}
Suppose from equation (\ref{e:fistequalll})  that $e_2(k_1)=0$. If from equation (\ref{e:secequalll}) we have that, in addition, $e_1(k_1)=0$, then 
$k_1$ is constant. On the other hand if from the same equation we conclude that $e_2(k_2)=0$, $S$ is rotationally symmetric. A similar argument holds for 
$e_1(k_2)=0$.
\end{pf}

\vspace{0.1in}

\begin{Note}
Proposition \ref{p:weingmargtra} can be extended to the case of the space $({\mathbb L}({\mathbb S}^3),{\mathbb G})$ of oriented geodesics in the three sphere, where ${\mathbb G}$ denotes the canonical neutral K\"ahler metric - further details can found in \cite{Alekseevsky} \cite{AA}.
\end{Note}

\vspace{0.1in}

\subsection{Surfaces of rank one}

Let $S\subset{\mathbb H}^3$ be an immersed surface with orthogonal geodesics $\Sigma\subset {\mathbb L}({\mathbb H}^3)$. 

\vspace{0.1in}

\begin{Prop}\label{p:lagrrankone}
If $\Sigma$ is of rank one then $S$ has constant principal curvature $\lambda=1$. 
\end{Prop}
\begin{pf}
See the Main Theorem of \cite{nikosbrendan}.
\end{pf}
\vspace{0.1in}

This Proposition, together with Proposition \ref{p:impaenri}, implies that every Lagrangian surface in ${\mathbb L}({\mathbb H}^3)$ of rank one is 
marginally trapped.

The following Theorem, extends this to all rank one surfaces:
\vspace{0.1in}
\begin{Thm}\label{t:allrankonemarg}
Every rank one surface in ${\mathbb L}({\mathbb H}^3)$ is marginally trapped.
\end{Thm}
\begin{pf}
Let $\Sigma$ be a  surface in ${\mathbb L}({\mathbb H}^3)$ of  rank one, locally parameterised by $\mu_1=\mu_1(s)$ and $\mu_2=\mu_2(s,t)$. The metric ${\mathbb G}$ induced on $\Sigma$ has the following components:
\[
g_{ss}=2{\mbox Im}\left[\frac{\partial_{s}\mu_1\partial_{s}\bar{\mu}_2}{(1+\mu_1\bar{\mu}_2)^2}\right]\qquad
g_{st}={\mbox Im}\left[\frac{\partial_{s}\mu_1\partial_{t}\bar{\mu}_2}{(1+\mu_1\bar{\mu}_2)^2}\right]\qquad
g_{tt}=0,
\]
which implies that $g$ is Lorentzian.
The only non-vanishing Christoffel symbols are:
\[
\Gamma^{s}_{ss}=g^{st}\partial_{s}g_{st}-\frac{1}{2}g^{st}\partial_{t}g_{ss},\qquad
\Gamma^{t}_{ss}=g^{tt}\partial_{s}g_{st}+\frac{1}{2}g^{st}\partial_{s}g_{ss}-\frac{1}{2}g^{tt}\partial_{t}g_{ss},
\]
\[
\Gamma^{t}_{st}=\frac{1}{2}g^{st}\partial_{t}g_{ss},\qquad\qquad \Gamma^{t}_{tt}=g^{st}\partial_{t}g_{st}.
\]
Observe that the induced metric $g$ of a rank one Lagrangian surface $\Sigma$ is scalar flat (see \cite{nikosbrendan}).

The second fundamental form $h=h_{ij}^{\mu_k}$ has non-vanishing components:  
\[
h_{ss}^{\mu_1}=\partial_{s}^2\mu_1-\frac{2\bar{\mu}_2(\partial_{s}\mu_1)^2}{1+\mu_1\bar{\mu}_2}-\Gamma_{ss}^{s}\partial_{s}\mu_1,
\quad 
h_{ss}^{\mu_2}=\partial_{s}^2\mu_2-\frac{2\bar{\mu}_1(\partial_{s}\mu_2)^2}{1+\bar{\mu}_1\mu_2}-\Gamma_{ss}^{s}\partial_{s}\mu_2-\Gamma_{ss}^{t}\partial_{t}\mu_2
\]
\[
h_{st}^{\mu_2}=\partial_{st}^2\mu_2-\frac{2\bar{\mu}_1\partial_{s}\mu_2\partial_{t}\mu_2}{1+\bar{\mu}_1\mu_2}-\Gamma_{st}^{t}\partial_{t}\mu_2,\qquad h_{tt}^{\mu_2}=\partial_{t}^2\mu_2-\frac{2\bar{\mu}_1(\partial_{t}\mu_2)^2}{1+\bar{\mu}_1\mu_2}-\Gamma_{tt}^{t}\partial_{t}\mu_2,
\]
with $h_{ij}^{\bar{\mu}_{k}}=\overline{h_{ij}^{\mu_{k}}}$.

The mean curvature vector $\vec{H}=2{\mathbb R}e(H^{\mu_1}\partial/\partial\mu_1+H^{\mu_2}\partial/\partial\mu_2)$  is:
\[
H^{\mu_{i}}=g^{ss}h^{\mu_{i}}_{ss}+2g^{st}h^{\mu_{i}}_{st}+g^{tt}h^{\mu_{i}}_{tt}.
\]
In particular,
\[
H^{\mu_1}=g^{ss}h^{\mu_1}_{ss}+2g^{st}h^{\mu_1}_{st}+g^{tt}h^{\mu_1}_{tt}=0,
\]
Thus
\[
\vec{H}=2{\mathbb R}e\Big(H^{\mu_2}\frac{\partial}{\partial\mu_2}\Big),
\]
so that ${\mathbb G}(\vec{H},\vec{H})=0$.
\end{pf}

\vspace{0.1in}

\subsection{Surfaces of rank two}

We now investigate marginally trapped Lagrangian surfaces of rank two. Let $\Sigma$ be a surface given by the graph 
$\mu_2=\mu_2(\mu_1,\bar\mu_1)$ and denote by $F$ the following expression $F=\bar\mu_2/(1+\mu_1\bar\mu_2)$. The Lagrangian condition is then 
\begin{equation}\label{e:lagrncon}
\overline\partial F=\partial\bar F,
\end{equation}
where $\partial$ and $\bar\partial$ denote differentiation with respect to $\mu_1$ and $\bar\mu_1$ respectively. 

\begin{Def}\label{d:potfunc}
The Lagrangian condition (\ref{e:lagrncon}) implies the existence of a function $h:\Sigma\rightarrow {\mathbb R}$ such that 
\begin{equation}\label{e:realfuncforlagr}
F=\partial h.
\end{equation}
The function $h$ is called the {\em potential function} of $\Sigma$. 

The distance along the normal geodesic between the point on the surface and the closest point to the origin in the ball model of ${\mathbb H}^3$ is 
a function $r:\Sigma\rightarrow {\mathbb R}$ called the {\em support function} of the Lagrangian surface $\Sigma$. 
\end{Def}

A brief computation shows the relationship between the support and potential functions is 
\[
r=\ln|F|+h+r_0=\ln|\partial h|+h+r_0,
\]
for some constant $r_0$.

\begin{Def}\label{d:lagangle}
For a Lagrangian graph $\mu_2=\mu_2(\mu_1,\bar\mu_1)$, define \emph{Lagrangian angle} $\phi$ of the surface by
\begin{equation}\label{e:tosigma00}
\frac{\partial\bar{\mu}_2}{(1+\mu_1\bar{\mu}_2)^2}=\sigma_0=|\sigma_0|e^{2i\phi}.
\end{equation}
\end{Def}

\vspace{0.1in}

\begin{Thm}\label{t:grapmargtra0}
Let $\Sigma$ be a rank two Lagrangian surface parameterised by $\mu_2=\mu_2(\mu_1,\bar\mu_1)$, with potential function $h$ and Lagrangian angle $\phi$.  
Then $\Sigma$ is marginally trapped iff $\partial e^{-i\phi+h}$ is either a real-valued  or an imaginary-valued function. 
\end{Thm}
\begin{pf}
The induced metric $g$ on $\Sigma$ is 
\[
g=-i\sigma_0d\mu_1^2+i\bar\sigma_0d\bar\mu_1^2,
\]
where, as before, $\sigma_0=\partial\bar\mu_2/(1+\mu_1\bar\mu_2)^2$. Non-degeneracy of the metric means that $\sigma_0\neq0$.

If $\vec{H}$ denotes the mean curvature vector, the vector field $J\vec{H}$ is 
tangential and is given by $J\vec{H}=a\partial+\bar a\overline\partial$, where 
\[
a=-\frac{1}{2\sigma_0}\left[\partial\ln\left(\frac{\bar\sigma_0}{\sigma_0}\right)-\frac{4\bar\mu_2}{1+\mu_1\bar\mu_2}\right].
\]
The surface is marginally trapped if and only if
\[
\partial e^{-i\phi}\pm\overline\partial e^{i\phi}=\mp e^{i\phi}\frac{\mu_2}{1+\bar\mu_1\mu_2}- e^{-i\phi}\frac{\bar\mu_2}{1+\mu_1\bar\mu_2},
\]
and since
\[
\partial h=F=\frac{\bar\mu_2}{1+\mu_1\bar\mu_2}
\]
this is equivalent to $\partial e^{-i\phi+h}=\mp\;\overline\partial e^{i\phi+h}$, proving the Theorem.
\end{pf}

\vspace{0.1in}

Note that the following holds
\begin{equation}\label{e:complpoinwithh}
\sigma_0=\partial^2 h+(\partial h)^2,
\end{equation}
and thus the relation between the potential function $h$ and the Lagrangian angle $\phi$ is 
\[
e^{4i\phi}=\frac{\partial^2 h+(\partial h)^2}{\overline\partial^2 h+(\overline\partial h)^2}.
\]

\vspace{0.1in}

The following Corollary gives examples of Lagrangian marginally trapped surfaces of rank two in ${\mathbb L}({\mathbb H}^3)$. 

\vspace{0.1in}

\begin{Cor}\label{t:grapmargtra}
If the potential function $h$ defined in  (\ref{e:realfuncforlagr}) is of the form 
\begin{equation}\label{e:thefunrealwithf}
h(\mu_1,\bar\mu_1)=f(\tau\mu_1+\bar\tau\bar\mu_1),
\end{equation}
where $f$ is a real function and $\tau\in{\mathbb S}^1$ is constant, then  the Lagrangian surface $\Sigma$ is marginally trapped.
\end{Cor}

\vspace{0.1in}

The next Proposition gives an equivalent condition of whether a marginally trapped surface Corollary's \ref{t:grapmargtra} form, is orthogonal to a 
Weingarten surface in ${\mathbb H}^3$.

\vspace{0.1in}

\begin{Prop}\label{p:ahgsu}
Let $\Sigma$ be the marginally trapped Lagrangian surface which is of the form given in Corollary {\ref{t:grapmargtra}}. Then $\Sigma$ is Weingarten, 
iff one of the following holds:
\begin{enumerate}
\item $\tau= 1$ or  $\tau=i$,
\item the real function $f$ defined by {(\ref{e:thefunrealwithf})} satisfies the following differential equation:
\[
4\mbox{Re}(\tau^2)[f''+(f')^2][f^{4}+2f'f^{3}+2(f'')^2]-(1+4\mbox{Re}(\tau^2))(f^{3}+2f'f'')^2=0.
\]
\end{enumerate}
Moreover, for any real constants $c_0\neq 0$ and $d_0$ the following support functions 
\begin{equation}\label{e:suporweingma1}
r(\mu_1,\bar\mu_1)=\ln\Big|\sinh[c_0(\tau\mu_1+\bar\tau\bar\mu_1)+d_0]\Big|+r_0,
\end{equation}
and 
\begin{equation}\label{e:suporweingma2}
r(\mu_1,\bar\mu_1)=\ln\Big|\sin[c_0(\tau\mu_1+\bar\tau\bar\mu_1)+d_0]\Big|+r_0,
\end{equation}
define local parameterisations of Weingarten surfaces in ${\mathbb H}^3$ such that the corresponding oriented geodesics are marginally 
trapped surfaces of rank two.
\end{Prop}
\begin{pf}
The fact that the induced metric is non-degenerate implies that $\sigma_0\neq 0$. Then by (\ref{e:complpoinwithh}) we have that 
\begin{equation}\label{e:consitiondik}
\partial^2h+(\partial h)^2\neq 0.
\end{equation}
Consider the real function $f$ satisfying (\ref{e:thefunrealwithf}). By (\ref{e:consitiondik}) we have $f''+(f')^2\neq 0$ and set 
$M=f''+(f')^2\neq 0$. The Gauss curvature of the Lagrangian surface $\Sigma$ given by the graph $\mu_2=\mu_2(\mu_1,\bar\mu_1)$ is \cite{gag}:
\begin{eqnarray}
K_{\Sigma}&=&\frac{i}{4|\sigma_0|^2} \Big(2(\bar\partial^2\sigma_0-\partial^2\bar\sigma_0)+\frac{(\partial\bar\sigma_0)^2-\bar\partial\bar\sigma_0\bar\partial\sigma_0}{\bar\sigma_0} -\frac{(\bar\partial\sigma_0)^2-\partial\sigma_0\partial\bar\sigma_0}{\sigma_0} \Big)\nonumber \\
&=&\frac{i(\tau^2-\bar\tau^2)}{4|\sigma_0|^2}\Big(2(\tau^2+\bar\tau^2)M''-(\tau^2+\bar\tau^2+1)\frac{(M')^2}{M}\Big).\nonumber
\end{eqnarray}
For a real non-zero constant $c_0$, the support functions given by (\ref{e:suporweingma1}) and (\ref{e:suporweingma2}) come from solving the equations $f''+(f')^2=c_0^2$ and $f''+(f')^2=-c_0^2$, respectively.
\end{pf}

\vspace{0.1in}

\begin{Note} 
We now use Proposition {\ref{p:ahgsu}} to construct examples of Lagrangian marginally trapped surfaces that are orthogonal to neither rotationally symmetric surfaces in ${\mathbb H}^3$ nor surfaces with constant principal curvature. For example, the marginally trapped surface, defined by 
\[
h(\mu_1,\bar\mu_1)=\left(\frac{1+i}{\sqrt{2}}\mu_1+\frac{1-i}{\sqrt{2}}\bar\mu_1\right)^2, 
\]
has Gauss curvature 
\[
K_{\Sigma}=\frac{4t^2}{(1+t^2)^3},
\]
where $t$ is the real-valued function:
\[
t=\frac{1+i}{\sqrt{2}}\mu_1+\frac{1-i}{\sqrt{2}}\bar\mu_1.
\]
The Lagrangian surface $\Sigma$ is orthogonal to the family of parallel surfaces $\{S_{r_0}\}_{r_0\in{\mathbb R}}$ in ${\mathbb H}^3$ given by the support function:
\[
r=\ln|\partial h|+h+r_0=\ln(\sqrt{2}|t|)+t^2+r_0.
\]
The fact that $\Sigma$ is not flat, implies that each surface $S_{r_0}$ is not Weingarten and thus is not rotationally symmetric. Moreover, the 
non-Weingarten condition also implies that none of the principal curvatures of $S_{r_0}$ are constant functions.
\end{Note}

\vspace{0.1in}

We are now in position to construct examples of Lagrangian spheres that are locally marginally trapped.
\begin{Prop}
Let $\Sigma$ be the Lagrangian graph $\mu_2=\mu_2(\mu_1,\bar\mu_1)$ and $h(\mu_1,\bar\mu_1)=f(\tau\mu_1+\bar\tau\bar\mu_1)$ be the real-valued function as 
defined in Corollary {\ref{t:grapmargtra}}. 
For every real constant $c\neq 0$ and integer $n\geq 2$, set $f'(x)=x/(x^{2n}+c^2)$, 
the corresponding Lagrangian marginally trapped surface $\Sigma_{n,c}$ is of topological type a sphere. 
\end{Prop}
\begin{pf}
We use polar coordinates $\mu_1=Re^{i\theta}$. The potential function $h$ can be described in these coordinates by
\[
h(R,\theta)=f(R\cos(\theta+\theta_0)),
\]
where $\theta_0$ is the constant with $\tau=e^{i\theta_0}\in {\mathbb S}^1$  and the real-valued function $f$ satisfies the following differential equation
\[
f'(x)=\frac{x}{x^{2n}+c^2}.
\] 
Without loss of generality, we assume that $\theta_0=0$. Then
\begin{equation}\label{e:immersionjjd}
\mu_2=\frac{f'}{1-Re^{-i\theta}f'},
\end{equation}
and we now prove that $\Sigma_{n,c}$ does not intersect with the reflected diagonal 
$\overline\Delta=\{(\mu_1,\mu_2)\in {\mathbb P}^1\times {\mathbb P}^1\; :\; 1+\mu_1\bar\mu_2=0\}$. From (\ref{e:immersionjjd}) we have
\begin{eqnarray}
1+\mu_1\bar\mu_2&=&(1-Re^{-i\theta}f')^{-1}\nonumber \\
&=& (1-Rf'\cos\theta-iRf'\sin\theta)^{-1}\label{e:contikj}
\end{eqnarray}
Set 
\[
a(R,\theta)=1-Rf'\cos\theta-iRf'\sin\theta,
\]
and then
\[
|a(R,\theta)|^2=\left(1-\frac{R^2\cos^2\theta}{c^2+R^{2n}\cos^{2n}\theta}\right)^2+\frac{R^4\sin^2 2\theta}{4(c^2+R^{2n}\cos^{2n}\theta)^2}.
\]
For $\theta\in [0,2\pi)\setminus\{\pi/2,3\pi/2\}$, we have
\[
|a(R,\pi/2)|^2=|a(R,3\pi/2)|^2=\lim_{R\rightarrow\infty}|a(R,\theta)|^2=\lim_{R\rightarrow 0}|a(R,\theta)|^2=1.
\]
This shows that $|a(R,\pi/2)|^2<\infty$ and therefore, from (\ref{e:contikj}) we conclude that $1+\mu_1\bar\mu_2\neq 0$. Thus, $\Sigma_{n,c}$ does not intersect the reflected diagonal and therefore it is compact. Moreover, from (\ref{e:immersionjjd}), the immersion of $\Sigma_{n,c}$ is given by
\[
\mu_2(R,\theta)=\frac{R\cos\theta}{c^2+R^{2n}\cos^{2n}\theta-R^{2}e^{-i\theta}\cos\theta}.
\]
Note that for $\theta_0\neq 0$, the immersion is
\[
\mu_2(R,\theta)=\frac{2Re^{-i\theta_0}\cos(\theta+\theta_0)}{c^2+2^{2n}R^{2n}\cos^{2n}(\theta+\theta_0)-2R^{2}e^{-i(\theta+\theta_0)}\cos(\theta+\theta_0)}.
\]
We have in any case, 
\[
\lim_{R\rightarrow\infty}\mu_2(R,\theta)=\lim_{R\rightarrow 0}\mu_2(R,\theta)=0,
\]
which shows that $\Sigma_{n,c}$ must be a sphere.
\end{pf}

\vspace{0.1in}

\begin{Ex}
For $n=2$ and any real constant $c\neq 0$, the orthogonal geodesics $\Sigma_{2,c}$ of the surfaces $S_{r_0}\subset {\mathbb H}^3$ that are parameterised 
by the support function 
\[
r=\ln\left[\frac{|\tau\mu_1+\bar\tau\bar\mu_1|}{(\tau\mu_1+\bar\tau\bar\mu_1)^4+c^2}\right]+{\textstyle{\frac{1}{2c}}}\tan^{-1}\left[\frac{(\tau\mu_1+\bar\tau\bar\mu_1)^2}{\sqrt{c}}\right]+r_0,
\]
are all marginally trapped Lagrangian spheres. For $n=3$ and $c=1$, the orthogonal geodesics of the surfaces $S'_{r_0}\subset {\mathbb H}^3$ that are 
parameterised by the support function 
\[
r=\ln\left[\frac{|\tau\mu_1+\bar\tau\bar\mu_1|}{(\tau\mu_1+\bar\tau\bar\mu_1)^6+1}\right]+{\textstyle{\frac{1}{6}}}\ln[1+(\tau\mu_1+\bar\tau\bar\mu_1)^2]
+{\textstyle{\frac{\sqrt{3}}{6}}}\tan^{-1}\left[\frac{2(\tau\mu_1+\bar\tau\bar\mu_1)^2-1}{\sqrt{3}}\right]
\]
\[
\qquad\qquad\qquad\qquad\qquad -{\textstyle{\frac{1}{12}}}\ln\left[1-(\tau\mu_1+\bar\tau\bar\mu_1)^2+(\tau\mu_1+\bar\tau\bar\mu_1)^4\right]+r_0,
\]
is a marginally trapped Lagrangian sphere. Here $S_{r_0}$ and $S'_{r_0}$  are not Weingarten surfaces, for any $\tau\neq 1,i$. Take for example 
$\tau=(1+i)/\sqrt{2}$ to see that the Gauss curvatures of $\Sigma_{2,c}$ and $\Sigma_{3,1}$ are not zero. Therefore, the surfaces 
$S_{r_0}$  and $S'_{r_0}$ in ${\mathbb H}^3$ orthogonal to $\Sigma_{2,c}$ and $\Sigma_{3,1}$ for $\tau=(1+i)/\sqrt{2}$ are not 
rotationally symmetric nor do they have a constant principal curvature. 
\end{Ex}

\vspace{0.3in}

\noindent {\bf Acknowledgements.} Nikos Georgiou is supported by Fapesp (2010/08669-9).

\vspace{0.3in}


\begin{thebibliography}{10}

\bibitem{Alekseevsky}
D.V. Alekseevsky, B. Guilfoyle and W. Klingenberg, {\it On the geometry of spaces of oriented geodesics} Ann. Global Anal. Geom.  {\bf 40}, 1--21 (2011).

\bibitem{AA} 
H. Anciaux, {\it Spaces of geodesics of pseudo-Riemannian space forms and normal congruences of hypersurfaces},  to appear in Transactions of the AMS .

\bibitem{ango} 
H. Anciaux and Y. Godoy {\it Marginally trapped submanifolds in Lorentzian space forms and in the Lorentzian product of the sphere by the real line}, (2012) [math.DG/1209.5118].

\bibitem{AGR} 
H. Anciaux, B. Guilfoyle, P. Romon, {\it Minimal submanifolds in the tangent bundle of a Riemannian surface},  J. Geometry and Physics. {\bf 61}, 237--247 (2011).

\bibitem{chendillen}
B.-Y. Chen, F. Dillen,{\it Classification of marginally trapped Lagrangian surfaces in Lorentzian complex forms}, J. Math. Phys. {\bf 48}, no. 1 (2007). 

\bibitem{gag}
N. Georgiou and B. Guilfoyle, {\it On the space of oriented geodesics of hyperbolic 3-space}, Rocky Mountain J. Math. {\bf 40}, 1183--1219 (2010).

\bibitem{nikosbrendan}
N. Georgiou and B. Guilfoyle, {\it A characterization of Weingarten surfaces in hyperbolic 3-space}, Abh. Math. Sem. Hamburg {\bf 80}, 233--253 (2010).

\bibitem{klbg}
B. Guilfoyle and W. Klingenberg, {\it An indefinite K\"ahler metric on the space of oriented lines}, J. London Math. Soc. {\bf 72}, 497--509 (2005).

\bibitem{salvai}
M. Salvai, {\it On the geometry of the space of oriented lines of hyperbolic space}, Glasg. Math. J.
{\bf 49}, 357--366 (2007).

\bibitem{palmer}
B. Palmer, {\it Area minimization among marginally trapped surfaces in Lorentz-Minkowski space}, Calc. Var. and Partial Diff. Eq. {\bf 41}, no. 3-4, 387--395, (2011).

\bibitem{penrose}
R. Penrose, {\it Gravitational collapse and space-time singularities}, Phys. Rev. Lett. {\bf 14}, 57--59 (1965).

\end{thebibliography}
\end{document}